\documentclass[10pt]{amsart}
\usepackage{amsmath, amssymb}
\usepackage{amsfonts}
\usepackage{mathrsfs}
\usepackage[color, arrow,matrix,curve,cmtip,ps]{xy}
\usepackage{paralist}
\usepackage{amsthm}
\usepackage{tikz-cd}
\usetikzlibrary{arrows,calc,matrix}
\tikzset{
curvarr/.style={
  to path={ -- ([xshift=2ex]\tikztostart.east)
    |- (#1) [near end]\tikztonodes
    -| ([xshift=-2ex]\tikztotarget.west)
    -- (\tikztotarget)}
  }
}

\usepackage[final]{pdfpages}
\usepackage{dsfont,txfonts}
\usepackage{tikz}
\usepackage{rotating}
\usepackage{mathrsfs}
\PassOptionsToPackage{hyphens}{url}\usepackage{hyperref}
\usepackage{enumitem}
\usepackage{graphicx}
\usepackage{mathtools}
\usetikzlibrary{arrows,chains,matrix,positioning,scopes}
\usepackage{listings}
\lstdefinelanguage{GAP}{%
  morekeywords={%
    Assert,Info,IsBound,QUIT,%
    TryNextMethod,Unbind,and,break,%
    continue,do,elif,%
    else,end,false,fi,for,%
    function,if,in,local,%
    mod,not,od,or,%
    quit,rec,repeat,return,%
    then,true,until,while%
  },%
  sensitive,%
  morecomment=[l]\#,%
  morestring=[b]",%
  morestring=[b]',%
}[keywords,comments,strings]

\usepackage[T1]{fontenc}
\usepackage[variablett]{lmodern}
\usepackage{xcolor}
\usepackage[citestyle=alphabetic,bibstyle=authortitle]{biblatex}
\usepackage{lipsum}

\newtheorem{theorem}{Theorem}[section]
\theoremstyle{definition}
\newtheorem{lemma}[theorem]{Lemma}

\newtheorem{remark}[theorem]{Remark}

\everymath{\displaystyle}
\allowdisplaybreaks

\newcommand{\set}[1]{\left\lbrace #1 \right\rbrace}

\begin{document}

\title{Computing the negative $K$-theory of finite groups of order $\leq 100$}
\author{Georg Lehner}
\subjclass[2020]{19D35}
\keywords{negative $K$-theory, group rings}
\begin{abstract}
We outline how the group $K_{-1}( \mathbb{Z}[G] )$ for a finite group $G$ can be computed using the computer language $GAP$ and compile a table of all groups $G$ of order less than $100$ that have torsion in $K_{-1}( \mathbb{Z}[G] )$ .
\end{abstract}
\maketitle

\section{Introduction}

Let $G$ be a finite group and let $\mathbb{Z}[G]$ be the integral group ring of $G$. The group
$$K_{-1} \mathbb{Z}[G]$$
naturally appears as a building block for the computation of higher $K$-theory and $A$-theory groups, for example via the Farrell-Jones conjecture \cite{Lueck2005SurveyOC}, and plays an important role in questions in geometric topology. An understanding of the groups $G$ for which $K_{-1} \mathbb{Z}[G]$ is non-trivial is thus useful both for calculations as well as for the construction of counterexamples. Fortunately, the negative $K$-theory of finite groups remains fairly accessible by means of representation theoretic tools, and a general description has been given by Carter.

\begin{theorem}[Carter, \cite{carternegktheory}] \label{negativektheory}
Let $G$ be finite. The group $K_{-1} \mathbb{Z} G$ has the form
$$K_{-1} \mathbb{Z} [G] = \mathbb{Z}^r \oplus (\mathbb{Z}/2)^s$$
where
$$ r = 1 - r_\mathbb{Q} + \sum_{p | ~ |G|}( r_{\mathbb{Q}_p} - r_{\mathbb{F}_p})$$
and $s$ is equal to the number of irreducible $\mathbb{Q}$-representations $I$ with even Schur index $m(I)$ but odd local Schur index $m_p(I)$ at every prime $p$ dividing the order of $G$.
\end{theorem}

Further computations of negative $K$-theory can be found in Lafont, Magurn, Ortiz \cite{LAFONT_2009} for the groups $D_n$, $D_n \times C_2$ and $A_5 \times C_2$, as well as in Magurn \cite{magurn_negative} for the groups $C_n$, $\text{Dic}_n$, $\tilde{T}$, $\tilde{O}$ and $\tilde{I}$. As is clear from Carter's theorem, a computation of the negative $K$-theory boils down to the computation of the values $r$ and $s$ for a concrete group $G$. We will show how to do this algorithmically and use the computer language $GAP$, \cite{GAP4}, to create tables of all groups of order $\leq 28$ together with the value of $K_{-1} \mathbb{Z} [G]$, as well as tables of all groups of order  $\leq 100$ for which $s > 0$.

\begin{remark}
A more conceptual understanding of the groups $K_{-1}\mathbb{Z}[G]$ as equivalence classes of \emph{singular characters} of $G$ is possible, and is discussed in previous work of the author \cite{lehner2021passage}.
\end{remark}

\subsection{Acknowledgements}

This work grew as part of the authors thesis. We want to thank Holger Reich, Wolfgang Lück, John Nicholson and Georg Loho for helpful commentary. Funded by the Deutsche Forschungsgemeinschaft (DFG, German Research Foundation) under Germany´s Excellence Strategy – The Berlin Mathematics Research Center MATH+ (EXC-2046/1, project ID: 390685689), and in particular the Berlin Mathematical School.

\section{Computing $K_{-1} \mathbb{Z} [G]$ for finite groups using GAP}


We want to give a description of how to compute the negative $K$-theory groups $K_{-1} \mathbb{Z} [G]$ using the computer algebra system $GAP$, \cite{GAP4}. The code used can be found online at \url{https://github.com/georglehner/negativektheory}.
In the following, fix a finite group $G$ and let $n$ be the order of $G$. As stated by theorem \ref{negativektheory}, we have
$$K_{-1} \mathbb{Z} [G] \cong \mathbb{Z}^{r(G)} \oplus \mathbb{Z}/2^{s(G)}$$
with the sum going over all prime numbers $p$ dividing the order of $G$. The two coefficients compute as
$$r(G) = 1 - r_\mathbb{Q} + \sum_{p | n}( r_{\mathbb{Q}_p} - r_{\mathbb{F}_p}) ,$$
where $r_k$ is defined as the number of irreducible representations of $G$ over the field $k$ and $s(G)$ is equal to the number of isomorphism classes of irreducible $\mathbb{Q}$-representations of $G$ with even global Schur index, but odd local Schur index at every prime $p$ dividing the order of $G$.

\subsection{The free rank $r(G)$}

To compute $r(G)$ we need to compute the number of irreducible representations of $G$ over the fields $\mathbb{Q}, \mathbb{F}_p$ and $\mathbb{Q}_p$. An essential tool for this is Berman's theorem. For the following, fix a field $k$ of characteristic $p$. We call an element $g$ of $G$ $p$-\emph{singular} if $p$ divides the order of $g$ and $p$-\emph{regular} otherwise. Note that if $p$ does not divide the order of $G$, all elements of $G$ are $p$-regular. Let $m$ be the exponent of $G$, i.e. the least common multiple over all orders of elements of $G$. Denote by $\zeta_m$ a primitive $m$-th root of unity. The Galois group $\text{Gal}(k(\zeta_m) : k)$ can by identified with a subgroup $T_m$ of $(\mathbb{Z}/m)^\times$, since any $\phi \in \text{Gal}(k(\zeta_m) : k)$ is uniquely defined by its image of $\zeta_m$, and so $\phi(\zeta_m) = \zeta_m^t$ for a unique $t \in T_m$. We say two elements $g_1$ and $g_2$ of $G$ are $k$-\emph{conjugate} if there exists $h \in G$ and a $t \in T_m$ such that $g_1^t = h g_2 h^{-1}$.

\begin{theorem}[Berman, \cite{berman}, see also \cite{reiner}]
The number of irreducible $k$-representations of $G$ is equal to the number of $k$-conjugacy classes of $p$-regular elements of $G$.
\end{theorem}

For the field $\mathbb{Q}$, the Galois group $\text{Gal}(\mathbb{Q}(\zeta_m) : \mathbb{Q})$ is always equal to the entire group $(\mathbb{Z}/m)^\times$. The number $r_\mathbb{Q}$ is thus equal to the number of conjugacy classes of cyclic subgroups of $G$.

In the case of the field $\mathbb{Q}_p$ a number $t$ is in  $t \in T_m$ iff $t$ is not divisible by $p$ and $t$ is congruent to a power of $p$ mod $\mu$, see \cite{greenberg2013local} Chapter IV, Section 4.

The following lemma will also be useful.
\begin{lemma}[Magurn \cite{magurn_negative}, see Lemma 1.] \label{Qpconjugacyclassessingular}
For each prime $p$ dividing $n$, $r_{\mathbb{Q}_p} - r_{\mathbb{F}_p}$ is equal to the number of $\mathbb{Q}_p$-conjugacy classes of $p$-singular elements in $G$.
\end{lemma}

\begin{remark} \label{pgroupstorsion} In particular, if the group $G$ is a $p$-group, all non-trivial elements of $G$ are $p$-singular, hence the number $r_{\mathbb{Q}_p} - r_{\mathbb{F}_p}$ is simply the number of $\mathbb{Q}_p$-conjugacy classes of non-trivial elements. In this case the exponent $m$ is a power of $p$, hence the Galois group $\text{Gal}(\mathbb{Q}_p(\zeta_m) : \mathbb{Q}_p)$ is actually the entire group $(\mathbb{Z}/m)^\times$. This means that $\mathbb{Q}$-conjugacy and $\mathbb{Q}_p$-conjugacy agree and thus
$$r(G) = 0. $$
More generally, $r(G) = 0$ for all groups $G$ that do not contain elements of non-prime power order.
\end{remark}

We will use the pre-existing functionality of $GAP$ to compute the set of conjugacy classes of $G$. For the following, we chose the function ConjugacyClassesByRandomSearch, but any of the predefined methods should work.

Since $k$-conjugacy is a coarser equivalence relation than conjugacy, it induces an equivalence relation on the set $cc$ of conjugacy classes of $G$. In order to simplify code, we have written a routine MyPartitionSet that takes a set \texttt{S} together with an equivalence relation \texttt{rel} as input and gives the set of equivalence classes as output.

\begin{lstlisting}[language=GAP]
MyPartitionSet := function( S, rel )
local Sprime, classes, s, allobjectsequivtos;
Sprime := S;
classes := [];
while not IsEmpty(Sprime) do
     s:= First(Sprime);
     allobjectsequivtos := Filtered ( Sprime , sprime -> rel(s,sprime) );
     Add(classes, allobjectsequivtos);
     SubtractSet(Sprime, allobjectsequivtos);
od;
return classes;
end;
\end{lstlisting}

We note that \texttt{rel} should be given as a function in two variables. The output is a list of \texttt{rel}-equivalence classes of elements of \texttt{S} .\footnote{We note that the code here presented only works if \texttt{rel} is actually an equivalence relation, and  {\texttt{MyPartitionSet}} does not check if this is satisfied.}

The set of $\mathbb{Q}$-conjugacy classes can then be computed with the following function \texttt{RationalConjugacyClasses}. This function takes a set \texttt{cc} of conjugacy classes as input and produces a list of sets of Q-conjugate conjugacy classes of G as output. We note that two elements of a group $G$ are rationally conjugate iff the cyclic subgroups they generate are conjugate.
\begin{lstlisting}[language=GAP]
RationalConjugacyClasses := function(cc)

local orders, m, galoisgroup, ccsortedbyorder, classesoforderi, kconjugate, 
rationalconjugacyclasses, c, t, allconjugatestoc;
orders := Set(List(cc, Representative), Order);
# The set of all appearing orders of elements of G
m := Lcm(orders);                             
# This is the exponent of G
galoisgroup := PrimeResidues(m);

rationalconjugacyclasses := [];
ccsortedbyorder := List( orders, i -> Filtered( cc , c -> Order(Representative(c)) = i ) );  
# ccSortedByOrder is a list with entries 
# the set of all conjugacy classes of given order i.
for classesoforderi in ccsortedbyorder do
    UniteSet(rationalconjugacyclasses, MyPartitionSet(classesoforderi, 
          {c,d} -> ForAny( galoisgroup , t -> Representative(c)^t in d ) ));
# We now partition all classes of order i into the set of equivalence classes
# of conjugacy classes with respect to k conjugacy
od;
return rationalconjugacyclasses;
end;
\end{lstlisting}

For the terms $r_{\mathbb{Q}_p} - r_{\mathbb{F}_p}$ we can simplify their computation using Lemma \ref{Qpconjugacyclassessingular} hence we need to compute the number of $\mathbb{Q}_p$-conjugacy classes of singular elements of $G$. Using Berman's theorem, the only thing left to understand is the Galois group $\text{Gal}(\mathbb{Q}_p(\zeta_m) : \mathbb{Q}_p)$. For this we need to compute $T_m$ as a subset of $(\mathbb{Z}/m)^\times$. We have that $t \in T_m$ iff $t$ is not divisible by $p$ and $t$ is congruent to a power of $p$ mod $\mu$.

In $GAP$ this can be realized with the following function, which gives the Galois group as a subset of prime residues of $m$.

\begin{lstlisting}[language=GAP]
GaloisGroupOfQpZetamOverQp := function(m, p) 

local q, mu, PowersOfPModmu, ZmUnits;
q := p^Length(Filtered(FactorsInt(m) , x -> x = p));
mu := BestQuoInt(m,q);
PowersOfPModmu := List( [1..OrderMod(p,mu)], i -> PowerModInt(p,i,mu) );
return Filtered( PrimeResidues(m), i -> i mod mu in PowersOfPModmu);
end;
\end{lstlisting}

Using this group we then define the equivalence relation of $\mathbb{Q}_p$-conjugacy on the set of conjugacy classes of $G$. Here \texttt{QpConjugate} is a function which takes as input a prime $p$, and two conjugacy classes \texttt{c} and \texttt{d} of $G$, and returns either \texttt{true} or \texttt{false}, depending on whether \texttt{c} and \texttt{d} are $\mathbb{Q}_p$-conjugate.

\begin{lstlisting}[language=GAP]
QpConjugate := function(p, c, d, m )
local Tm;
Tm := GaloisGroupOfQpZetamOverQp(m,p);        
return ForAny( Tm, t->Representative(c)^t in d) ;      
\end{lstlisting}

With this defined, we now compute the set of singular $\mathbb{Q}_p$-conjugacy classes. Recall that an element $g$ is $p$-singular if $p$ divides the order of $g$. The function \texttt{SingularQpConjugacyClasses} takes a set of conjugacyclasses \texttt{cc} and a prime $p$ as input, and produces a list of sets of $\mathbb{Q}_p$-conjugate conjugacy classes of $G$ with $p$-singular orders.

\begin{lstlisting}[language=GAP]
SingularQpConjugacyClasses := function(cc,p)

local orders, singularorders, m, ccsortedbyorder, classesoforderi, 
singularqpconjugacyclasses;
orders := Set(List(cc, c->Order(Representative(c)) )); 
# The set of all appearing orders of elements of G
m := Lcm(orders);                             
# This is the exponent of G

singularorders := Filtered(List(cc, c->Order(Representative(c)) ), i -> i mod p = 0 ); 
# We will only look at orders that are multiples of p

singularqpconjugacyclasses := [];
ccsortedbyorder := List( singularorders, 
     i -> Filtered( cc , c -> Order(Representative(c)) = i ) ); 
# ccSortedByOrder is a list with entries the set of
# all conjugacy classes of given order i. 
# We only care about orders divisible by p, hence SingularOrders
for classesoforderi in ccsortedbyorder do
    UniteSet(singularqpconjugacyclasses, MyPartitionSet(classesoforderi, 
         {c,d} -> QpConjugate(p,c,d,m) ) );
# We now partition all classes of order i into the set of equivalence classes
# of conjugacy classes with respect to Qp conjugacy
od;
return singularqpconjugacyclasses;
end;
\end{lstlisting}

Using the formula
$$r(G) = 1 - r_\mathbb{Q} + \sum_{p | n} \#\set{\text{singular } \mathbb{Q}_p\text{-conjugacy classes}}, $$
we can now compute $r(G)$ using the function:

\begin{lstlisting}[language=GAP]
rOfGroup := function(G)
local primes, rQ, rsingQp, cc;
cc := ConjugacyClassesByRandomSearch(G);
primes := PrimeDivisors(Order(G));
rQ := Length( RationalConjugacyClasses(cc) );
rsingQp := List( primes, p -> Length(SingularQpConjugacyClasses(cc,p)) );
return 1-rQ + Sum(rsingQp);
end;
\end{lstlisting}

\subsection{The $\mathbb{Z}/2$ rank $s(G)$}

The number $s(G)$ is computed as the number of irreducible rational representations of $G$ that have global even Schur index but odd Schur index at every prime $p$ dividing the order of $G$. This can be computed using the Wedderga package in $GAP$ \cite{Wedderga}. The main function we will use is WedderburnDecompositionInfo, which computes the Wedderburn decomposition of $G$ as a list in which each entry corresponds to a numerical description of the cyclotomic algebras appearing in the Wedderburn decomposition. The global Schur index can be computed with the function SchurIndex and the local Schur indices with LocalIndicesOfCyclotomicAlgebra, which gives a list of pairs $(p,m_p)$ contain to the finitely many places $p$ at which the local Schur index $m_p$ is non-trivial. Note that $p = \infty$ corresponding to the real Schur index can also appear in this list, which we have to remove when computing the relevant indices for $s(G)$. The following function \texttt{AlgebraContributesToS} checks if a cyclotomic algebra $A$ described in the format given by WedderburnDecompositionInfo gives a contribution to the number $s(G)$, and returns a truth value. We note that $A$ must be in the form of an output of the function
\texttt{WedderburnDecompositionInfo} of the package \texttt{wedderga}, and \texttt{primes} is a list of primes that the function checks against.

\begin{lstlisting}[language=GAP]

AlgebraContributesToS := function(A,primes)
local RelevantIndices;
if IsEvenInt(SchurIndex(A)) then
RelevantIndices := Filtered( LocalIndicesOfCyclotomicAlgebra(A) , l -> l[1] in primes );
return ForAll(RelevantIndices, l -> IsOddInt( l[2] ) );
else return false; fi;
end;
\end{lstlisting}

With AlgebraContributesToS, it is now a simple matter of counting all the summands in the Wedderburn decomposition which contribute.

\begin{lstlisting}[language=GAP]
sOfGroup := function(G)
local n, primes, ww, s;
n := Order(G);
primes := PrimeDivisors(n);
ww := WedderburnDecompositionInfo( GroupRing( Rationals, G ) );
s := Number( ww, A -> AlgebraContributesToS(A,primes) );
return s;
end;
\end{lstlisting}

\subsection{Wrapping things up}

For user-friendliness we can put both rOfGroup and sOfGroup into the simple function KMinusOne which returns the values of r and s in one record.

\begin{lstlisting}[language=GAP]
KMinusOne := function(G)
return rec(r := rOfGroup(G), s := sOfGroup(G) );
end;
\end{lstlisting}

As an example, take the binary icosahedral group $\tilde{I}$, also identifiable as $SL(2,\mathbb{F}_5)$. The ID of this group in the small groups library is given as $(120,5)$.

\begin{lstlisting}[language=GAP]
gap> G := SmallGroup(120,5); 
Group([ (1,2,4,8)(3,6,9,5)(7,12,13,17)(10,14,11,15)(16,20,21,24)(18,22,19,23), 
 (1,3,7)(2,5,10)(4,9,13)(6,11,8), (12,16,20)(14,18,22)(15,19,23)(17,21,24) ])
gap> StructureDescription(G);
"SL(2,5)"
gap> KMinusOne(G);
rec( r := 2, s := 1 )
\end{lstlisting}

From this result we can read off that $K_{-1} \mathbb{Z}[SL(2,\mathbb{F}_5)] = \mathbb{Z}^2 \oplus
 \mathbb{Z}/2$, which agrees with the results of \cite{magurn_negative}.

\subsection{$K_{-1} \mathbb{Z} [G]$ of all groups of order $\leq 28$} \label{computingktheory}
 
In the following we provide tables of all groups of order $\leq 28$ together with the value of their negative $K$-theory groups $K_{-1} \mathbb{Z} [G]$. In order to compute a table of all groups with non-trivial $s$, we use the function KMinusOne constructed in the previous section as well as the functionality provided by the SmallGroup library in $GAP$. The SmallGroups library allows one to address all groups of reasonably small order by two parameters, the order of the group $n$ as well as an index $i$. The idea is to fix an order $n$, and then compute $s$ for each group of order $n$ by going through all possible indices $i$. In the column ``Structure'' we describe the group by a standard name or, in case there is no such name, we include a printout of the functionality StructureDescription of $GAP$. Note that StructureDescription(G) does not specify the group $G$ up to isomorphism, but only serves to give a quick idea of the type of group which $G$ is. The relevant tables are Tables \ref{groupsorderless12}, \ref{groupsorderless24} and \ref{groupsorder24}.
 
\begin{remark}
We want to summarize two conditions for the vanishing of negative $K$-theory for the convenience of the reader.
\begin{itemize}
\item The free rank $r$ vanishes iff $G$ contains no element of non-prime power order. This is visible for example in the case of $n = 6$, where the group $C_6$ has $K_{-1} \mathbb{Z}[C_6] = \mathbb{Z}$, but for $S_3$ we have $K_{-1} \mathbb{Z}[S_3] = 0.$
\item The torsion rank $s$ vanishes if the group order is not a multiple of $4$.\footnote{This is because any irreducible $\mathbb{Q}$-representation that contributes to $s$ splits into non-trivial quaternionic representations of $G$ over $\mathbb{R}$, and the existence of these implies that the group order is a multiple of $4$, see \cite{Unger_2019} Theorem 4.1. in conjunction with Theorem 2.2, (2)}   This is however not a sufficient criterion. For example, there do not exist groups with $s > 0$ of order 4, 8, 12, 28, 36, 44, 76 or 92. Purely experimentally, among the groups with order a multiple of $4$, the groups with $s > 0$ tend to be rare, although that could be down to the fact that we are only looking at small values of $n$.
\end{itemize}
\end{remark}

\begin{table} 
\begin{center}
\begin{tabular}{| c | c | c | c | c | c | c |}
\hline
$n = 1$ & Index  &  Structure  & $r(G)$ & $s(G)$ & $K_{-1}\mathbb{Z}[G]$  \\ 
\hline
 & $1$ & $1$ & 0 & 0 &   $0$  \\
\hline
\hline
$n = 2$ & Index  &  Structure  & $r(G)$ & $s(G)$ & $K_{-1}\mathbb{Z}[G]$  \\ 
\hline
 & $1$ & $C_2$ & 0 & 0 &  $0$   \\
\hline
\hline
$n = 3$ & Index  &  Structure  & $r(G)$ & $s(G)$ & $K_{-1}\mathbb{Z}[G]$  \\ 
\hline
 & $1$ & $C_3$ & 0 & 0 &  $0$   \\
\hline
\hline
$n = 4$ & Index  &  Structure  & $r(G)$ & $s(G)$ & $K_{-1}\mathbb{Z}[G]$  \\ 
\hline
 & $1$ & $C_4$ & 0 & 0 &  $0$   \\
 & $2$ & $C_2 \times C_2$ & 0 & 0 &  $0$   \\
\hline
\hline
$n = 5$ & Index  &  Structure  & $r(G)$ & $s(G)$ & $K_{-1}\mathbb{Z}[G]$  \\ 
\hline
 & $1$ & $C_5$ & 0 & 0 &  $0$   \\
\hline
\hline
$n = 6$ & Index  &  Structure  & $r(G)$ & $s(G)$ & $K_{-1}\mathbb{Z}[G]$  \\ 
\hline
 & $1$ & $S_3$ & 0 & 0 &  $0$   \\
 & $2$ & $C_6$ & 1 & 0 &  $\mathbb{Z}$   \\
\hline
\hline
$n = 7$ & Index  &  Structure  & $r(G)$ & $s(G)$ & $K_{-1}\mathbb{Z}[G]$  \\ 
\hline
 & $1$ & $C_7$ & 0 & 0 &  $0$   \\
\hline
\hline
$n = 8$ & Index  &  Structure  & $r(G)$ & $s(G)$ & $K_{-1}\mathbb{Z}[G]$  \\ 
\hline
 & $1$ & $C_8$ & 0 & 0 &  $0$   \\
 & $2$ & $C_4 \times C_2$ & 0 & 0 &  $0$   \\
 & $3$ & $D_8$ & 0 & 0 &   $0$  \\
 & $4$ & $Q_8$ & 0 & 0 &  $0$   \\
 & $5$ & $C_2 \times C_2 \times C_2$ & 0 & 0 &  $0$   \\
\hline
\hline
$n = 9$ & Index  &  Structure  & $r(G)$ & $s(G)$ & $K_{-1}\mathbb{Z}[G]$  \\ 
\hline
 & $1$ & $C_9$ & 0 & 0 &   $0$  \\
 & $2$ & $C_3 \times C_3$ & 0 & 0 &  $0$   \\
\hline
\hline
$n = 10$ & Index  &  Structure  & $r(G)$ & $s(G)$ & $K_{-1}\mathbb{Z}[G]$  \\ 
\hline
 & $1$ & $D_{10}$ & 0 & 0 &  $0$   \\
 & $2$ & $C_{10}$ & 1 & 0 &  $\mathbb{Z}$   \\
\hline
\hline
$n = 11$ & Index  &  Structure  & $r(G)$ & $s(G)$ & $K_{-1}\mathbb{Z}[G]$  \\ 
\hline
 & $1$ & $C_{11}$ & 0 & 0 &  $0$   \\
\hline
\hline
$n = 12$ & Index  &  Structure  & $r(G)$ & $s(G)$ & $K_{-1}\mathbb{Z}[G]$  \\ 
\hline
 & $1$ & $\mathrm{Dic}_{3}$ & 1 & 0 &  $\mathbb{Z}$   \\
 & $2$ & $C_{12}$ & 2 & 0 &  $\mathbb{Z}^2$   \\
 & $3$ & $A_{4}$ & 0 & 0 &  $0$   \\
 & $4$ & $D_{12}$ & 1 & 0 &  $\mathbb{Z}$   \\
 & $5$ & $C_{6} \times C_2$ & 3 & 0 &  $\mathbb{Z}^3$   \\
\hline
\hline
$n = 13$ & Index  &  Structure  & $r(G)$ & $s(G)$ & $K_{-1}\mathbb{Z}[G]$  \\ 
\hline
 & $1$ & $C_{13}$ & 0 & 0 &  $0$   \\
\hline
\hline
$n = 14$ & Index  &  Structure  & $r(G)$ & $s(G)$ & $K_{-1}\mathbb{Z}[G]$  \\ 
\hline
 & $1$ & $D_{14}$ & 0 & 0 &  $0$   \\
 & $1$ & $C_{14}$ & 2 & 0 &  $\mathbb{Z}^2$    \\
\hline
\hline
$n = 15$ & Index  &  Structure  & $r(G)$ & $s(G)$ & $K_{-1}\mathbb{Z}[G]$  \\ 
\hline
 & $1$ & $C_{15}$ & 0 & 0 &  $0$   \\
\hline
\end{tabular}  
\caption{All groups $G$ of order $ n < 16$ and the value of $K_{-1}\mathbb{Z}[G]$.} \label{groupsorderless12}
\end{center}
\end{table}

\begin{table} 
\begin{center}
\begin{tabular}{| c | c | c | c | c | c | c |}
\hline
$n = 16$ & Index  &  Structure  & $r(G)$ & $s(G)$ & $K_{-1}\mathbb{Z}[G]$  \\ 
\hline
 & $1$ & $C_{16}$ & 0 & 0 &   $0$  \\
 & $2$ & $C_4 \times C_4$ & 0 & 0 &   $0$  \\
 & $3$ & (C4 x C2) : C2 & 0 & 0 &   $0$  \\
 & $4$ & C4 : C4 & 0 & 0 &   $0$  \\
 & $5$ & $C_8 \times C_2$ & 0 & 0 &   $0$  \\
 & $6$ & C8 : C2 & 0 & 0 &   $0$  \\
 & $7$ & $D_{16}$ & 0 & 0 &   $0$  \\
 & $8$ & $QD_{16}$ & 0 & 0 &   $0$  \\
 & $9$ & $Q_{16}$ & 0 & 1 &   $\mathbb{Z}/2$  \\
 & $10$ & $C_4 \times C_2 \times C_2$ & 0 & 0 &   $0$  \\
 & $11$ & $C_2 \times D_8$ & 0 & 0 &   $0$  \\
 & $12$ & $C_2 \times Q_8$ & 0 & 0 &   $0$  \\
 & $13$ & (C4 x C2) : C2, Pauli group & 0 & 0 &   $0$  \\
 & $14$ & $C_2 \times C_2 \times C_2 \times C_2$ & 0 & 0 &   $0$  \\
\hline
\hline
$n = 17$ & Index  &  Structure  & $r(G)$ & $s(G)$ & $K_{-1}\mathbb{Z}[G]$  \\ 
\hline
 & $C_{17}$ & $1$ & 0 & 0 &   $0$  \\
\hline
\hline
$n = 18$ & Index  &  Structure  & $r(G)$ & $s(G)$ & $K_{-1}\mathbb{Z}[G]$  \\ 
\hline
 & $1$ & $D_{18}$ & 0 & 0 &   $0$  \\
 & $2$ & $C_{18}$ & 2 & 0 &   $\mathbb{Z}^2$  \\
 & $3$ & $C_3 \times S_3$ & 1 & 0 &   $\mathbb{Z}$  \\
 & $4$ & (C3 x C3) : C2 & 0 & 0 &   $0$  \\
 & $5$ & $C_6 \times C_3$ & 4 & 0 &   $\mathbb{Z}^4$  \\
\hline
\hline
$n = 19$ & Index  &  Structure  & $r(G)$ & $s(G)$ & $K_{-1}\mathbb{Z}[G]$  \\ 
\hline
 & $1$ & $C_{19}$ & 0 & 0 &   $0$  \\
\hline
\hline
$n = 20$ & Index  &  Structure  & $r(G)$ & $s(G)$ & $K_{-1}\mathbb{Z}[G]$  \\ 
\hline
 & $1$ & $\mathrm{Dic}_5$ & 1 & 1 &   $\mathbb{Z} \oplus \mathbb{Z}/2$  \\
 & $2$ & $C_{20}$ & 3 & 0 &   $\mathbb{Z}^3$  \\
 & $3$ & $F_5$ & 0 & 0 &   $0$  \\
 & $4$ & $D_{20}$ & 1 & 0 &   $\mathbb{Z}$  \\
 & $5$ & $C_{10} \times C_2$ & 3 & 0 &   $\mathbb{Z}^3$  \\
\hline
\hline
$n = 21$ & Index  &  Structure  & $r(G)$ & $s(G)$ & $K_{-1}\mathbb{Z}[G]$  \\ 
\hline
 & $1$ & C7 : C3 & 0 & 0 &   $0$  \\
 & $2$ & $C_{21}$ & 2 & 0 &   $\mathbb{Z}^2$  \\
\hline
\hline
$n = 22$ & Index  &  Structure  & $r(G)$ & $s(G)$ & $K_{-1}\mathbb{Z}[G]$  \\ 
\hline
 & $1$ & $D_{22}$ & 0 & 0 &   $0$  \\
 & $2$ & $C_{22}$ & 1 & 0 &   $\mathbb{Z}$  \\
\hline
\hline
$n = 23$ & Index  &  Structure  & $r(G)$ & $s(G)$ & $K_{-1}\mathbb{Z}[G]$  \\ 
\hline
 & $1$ & $C_{23}$ & 0 & 0 &   $0$  \\
\hline
\end{tabular}  
\caption{All groups $G$ of order $ 16 \leq n <24$ and the value of $K_{-1}\mathbb{Z}[G]$.} \label{groupsorderless24}
\end{center}
\end{table}

\begin{table} 
\begin{center}
\begin{tabular}{| c | c | c | c | c | c | c |}
\hline
$n = 24$ & Index  &  Structure  & $r(G)$ & $s(G)$ & $K_{-1}\mathbb{Z}[G]$  \\ 
\hline
 & $1$ & C3 : C8 & 2 & 0 &   $\mathbb{Z}^2$  \\
 & $2$ & $C_{24}$ & 4 & 0 &   $\mathbb{Z}^4$  \\
 & $3$ & $SL(2,3)$ & 1 & 0 &   $\mathbb{Z}$  \\
 & $4$ & C3 : Q8 & 2 & 1 &   $\mathbb{Z}^2 \oplus \mathbb{Z}/2$  \\
 & $5$ & $C_4 \times S_3$ & 2 & 0 &   $\mathbb{Z}^2$  \\
 & $6$ & $D_{24}$ & 2 & 0 &   $\mathbb{Z}^2$  \\
 & $7$ & C2 x (C3 : C4) & 3 & 0 &   $\mathbb{Z}^3$  \\
 & $8$ & (C6 x C2) : C2 & 2 & 0 &   $\mathbb{Z}^2$  \\
 & $9$ & $C_{12} \times C_2$ & 5 & 0 &   $\mathbb{Z}^5$  \\
 & $10$ & $C_3 \times D_8$ & 4 & 0 &   $\mathbb{Z}^4$  \\
 & $11$ & $C_3 \times Q_8$ & 4 & 0 &   $\mathbb{Z}^4$  \\
 & $12$ & $S_{4}$ & 0 & 0 &   $0$  \\
 & $13$ & $C_{2} \times A_4$ & 1 & 0 &   $\mathbb{Z}$  \\
 & $14$ & $C_2 \times C_2 \times S_3$ & 3 & 0 &   $\mathbb{Z}^3$  \\
 & $15$ & $C_6 \times C_2 \times C_2$ & 7 & 0 &   $\mathbb{Z}^7$  \\
\hline
\hline
$n = 25$ & Index  &  Structure  & $r(G)$ & $s(G)$ & $K_{-1}\mathbb{Z}[G]$  \\ 
\hline
 & $1$ & $C_{25}$ & 0 & 0 &   $0$  \\
 & $2$ & $C_{5} \times C_{5}$ & 0 & 0 &   $0$  \\
\hline
\hline
$n = 26$ & Index  &  Structure  & $r(G)$ & $s(G)$ & $K_{-1}\mathbb{Z}[G]$  \\ 
\hline
 & $1$ & $D_{26}$ & 0 & 0 &   $0$  \\
 & $2$ & $C_{26}$ & 1 & 0 &   $\mathbb{Z}$  \\
\hline
\hline
$n = 27$ & Index  &  Structure  & $r(G)$ & $s(G)$ & $K_{-1}\mathbb{Z}[G]$  \\ 
\hline
 & $1$ & $C_{27}$ & 0 & 0 &   $0$  \\
 & $2$ & $C_9 \times C_3$ & 0 & 0 &   $0$  \\
 & $3$ & (C3 x C3) : C3, Heisenberg group & 0 & 0 &   $0$  \\
 & $4$ & C9 : C3, Extraspecial group  & 0 & 0&   $0$  \\
 & $5$ & $C_3 \times C_3 \times C_3$ & 0 & 0 &   $0$  \\
\hline
\hline
$n = 28$ & Index  &  Structure  & $r(G)$ & $s(G)$ & $K_{-1}\mathbb{Z}[G]$  \\ 
\hline
 & $1$ & $\mathrm{Dic}_7$ & 1 & 0 &   $\mathbb{Z}$  \\
 & $2$ & $C_{28}$ & 4 & 0 &   $\mathbb{Z}^4$  \\
 & $3$ & $D_{28}$ & 1 & 0 &   $\mathbb{Z}$  \\
 & $4$ & $C_{14} \times C_2$ & 6 & 0 &   $\mathbb{Z}^6$  \\
\hline
\end{tabular}  
\caption{All groups $G$ of order $24 \leq n \leq 28$ and the value of $K_{-1}\mathbb{Z}[G]$.} \label{groupsorder24}
\end{center}
\end{table} 

\clearpage

\subsection{All groups of order $\leq 100$ with non-trivial torsion in $K_{-1} \mathbb{Z} [G]$} \label{computings}

We will also include a computation of all finite groups $G$ with $s(G) > 0 $ of order $n$ less than $100$. The relevant tables are Tables \ref{groupsorderless64}, \ref{groupsorder64}, \ref{groupsorder65to95} and \ref{groupsorder96to100}. However, a few comments are needed. First of all, torsion can only appear in $K_{-1} Z[G]$ when the order $n$ of the group $G$ is divisible by $4$. Moreover, we want to highlight that the binary polyhedral groups
$$ \text{Dic}_{n}, \tilde{O}, \tilde{I}$$
play a special role as examples of groups with $s > 0$. These are finite subgroups of $SU(2)$, called dicyclic group of order $n$ for $n \geq 4$, binary octahedral and binary icosahedral group, respectively. They are defined as the inverse images under the double covering $SU(2) \rightarrow SO(3)$ of the groups
$$ C_{2n}, O, I$$
with $C_{2n}$ the cyclic group of order $2n$, $O$ being the octahedral group and $I$ the icosahedral group, thought of as the symmetries of the regular $2n$-gon, cube and icosahedron, respectively. Magurn \cite{magurn_negative} computes the negative $K$-theory groups of all of these. The smallest group with $s>0$ is the group $Q_{16} = \text{Dic}_4$. Note that if $n = 2^k$ is a power of $2$, the group $\text{Dic}_{n}$ is also called generalized quaternion group $Q_{2^{k+2}}$ of order $4 \cdot 2^k = 2^{k+2}$.

Similarly to the previous section, in order to compute a table of all groups with non-trivial $s$, we use the function sOfGroup constructed in the previous section as well as the functionality provided by the SmallGroup library in $GAP$. We will only print groups with non-trivial $s$. If the group $G$ is of the form $ \text{Dic}_{n}, \tilde{O}$, $\tilde{I}$, or a product of these and another group, we address the group as such in the column ``Structure''. If it is not of this form, we include a printout of the functionality StructureDescription of $GAP$. Just as in the previous section, we want to warn the reader that StructureDescription(G) does not specify the group $G$ up to isomorphism, but only serves to give a quick idea of the type of group which $G$ is.

If $G$ has a quotient $G/N$ with $s(G/N) > 0$, then each irreducible $G/N$-representation $(I, \rho)$ that contributes to $s(G/N)$ gives a $G$-representation $(I, \rho \circ \pi)$ with $\pi  \colon G \rightarrow G/N$, which contributes to $s$. In the column ``Quotients'' we have included all smaller groups with $s > 0$ which appear as quotients. Most groups which are not binary polyhedral groups which appear in our list are, in fact, explained by having binary polyhedral groups as quotients, i.e. the group ring $\mathbb{Z} [G]$ fails the Eichler condition.\footnote{The Eichler condition on the group ring $\mathbb{Z} [G]$ is the condition that $ \mathbb{Z} [G] \otimes \mathbb{R} \cong \mathbb{R}G$ contains no copies of the quaternions $\mathbb{H}$ in its Wedderburn decomposition, see \cite{swan1970k}, Chapter 9., also \cite{NICHOLSON_2020}.}

However, not all groups in our list have their irreducible representations contributing to $s$ lifted from binary polyhedral quotients.  The smallest example is the group $C_8.C^2_2$ of order $32$ with SmallGroup library ID $(32,44)$. It has a unique irreducible representation contributing to $s$ induced from the subgroup $Q_{16}$. We included in Tables \ref{groupswithnoquotients128} and \ref{groupswithnoquotients180} all groups of order $\leq 180$ which do not have a quotient with non-trivial $s$.

\begin{table} 
\begin{center}
\begin{tabular}{| c | c | c | c | c | c | }
\hline
$n = 16$ & Index  &  Structure  & $s(G)$  & Quotients \\ 
\hline
 & 9 & $Q_{16} = \text{Dic}_4$  & 1 & \\
\hline
\hline
$n = 20$ & Index  &  Structure  & $s(G)$  & Quotients \\ 
\hline
 & 1 & $\text{Dic}_5$ & 1 &  \\
\hline
\hline
$n = 24$ & Index  &  Structure  & $s(G)$ & Quotients \\ 
\hline
 & 4 & $\text{Dic}_6$ & 1  \\
\hline
\hline
$n = 28$ & Index  &  Structure  & $s(G)$ & Quotients  \\ 
\hline
\hline
$n = 32$ & Index  &  Structure  & $s(G)$ & Quotients  \\ 
\hline
 & 10 & Q8 : C4 & 1 & $Q_{16}$  \\
 & 14 & C8 : C4 & 1 & $Q_{16}$  \\
 & 20 & $Q_{32} = \text{Dic}_8$ & 1 &  \\
 & 41 & $C_2 \times Q_{16}$ & 2 &  \\
 & 44 & (C2 x Q8) : C2 = $C_8.C_2^2$ & 1 & \\
\hline
 \hline
$n = 36$ & Index  &  Structure  & $s(G)$ & Quotients  \\ 
\hline
\hline
$n = 40$ & Index  &  Structure  & $s(G)$ & Quotients  \\ 
\hline
 & 1 & C5 : C8 & 1 & $\text{Dic}_5$  \\
 & 4 & $\text{Dic}_{10}$ & 1 & \\
 & 7 & $C_2 \times \text{Dic}_5$ & 2 & \\
\hline
\hline
$n = 44$ & Index  &  Structure  & $s(G)$ & Quotients  \\ 
\hline
\hline
$n = 48$ & Index  &  Structure  & $s(G)$ & Quotients  \\ 
\hline
 & 8 & $\text{Dic}_{12}$ & 2 & \\
 & 12 & (C3 : C4) : C4 & 1 & $\text{Dic}_{6}$ \\
 & 13 & C12 : C4 & 1 & $\text{Dic}_{6}$ \\
 & 18 & C3 : Q16 & 1 & $Q_{16}$ \\
 & 27 & $C_3 \times Q_{16}$ & 1 & \\
 & 28 & $\tilde{O}$ & 1 & \\
 & 34 & $C_2 \times \text{Dic}_6$ & 2 & \\
\hline
\hline
$n = 52$ & Index  &  Structure  & $s(G)$ & Quotients \\ 
\hline
 & 1 & $\text{Dic}_{13}$& 1 & \\
\hline
\hline
$n = 56$ & Index  &  Structure  & $s(G)$ & Quotients \\ 
\hline
 & 3 & $\text{Dic}_{14}$ & 1 & \\
\hline
\hline
$n = 60$ & Index  &  Structure  & $s(G)$ & Quotients  \\ 
\hline
 & 2 & $C_3 \times \text{Dic}_{5}$ & 1 & \\
 & 3 & $\text{Dic}_{15}$ & 2 & \\
\hline
\end{tabular}  
\caption{All groups $G$ with $s(G)>0$ of order $ n < 64$.} \label{groupsorderless64}
\end{center}
\end{table}

\begin{table}
\begin{center}
\begin{tabular}{| c | c | c | c | c | }
\hline
n = 64 & Index  &  Structure  & s(G) & Quotients \\ 
\hline
 & 7 & Q8 : C8 & 1 & $Q_{16}$ \\
 & 9 & (C4 : C4) : C4 & 1 & $Q_{16}$ \\
 & 13 & (C2 x C2) . ((C4 x C2) : C2) & 1 & $Q_{16}$ \\
 & 14 & (C2 x C2) . ((C4 x C2) : C2) & 2 & 2 copies of $Q_{16}$ \\
 & 16 & C8 : C8 & 1 & $Q_{16}$ \\
 & 21 & (C4 : C4) : C4 & 1 & $Q_{16}$ \\
 & 39 & Q16 : C4 & 1 & $Q_{32}$ \\
 & 43 & C2 . ((C8 x C2) : C2) & 1 & \\
 & 47 & C16 : C4 & 2 & $Q_{16}, Q_{32}$ \\
 & 48 & C16 : C4 & 1 &  $Q_{16}$ \\
 & 49 & C4 . D16 = C8 . (C4 x C2) & 1 & $Q_{16}$ \\
 & 54 & $Q_{64} = \text{Dic}_{16}$ & 1 & \\
 & 96 & C2 x (Q8 : C4) & 2 & \\
 & 107 & C2 x (C8 : C4) & 2 & \\
 & 120 & C4 x Q16 & 2 & \\
 & 129 & (C2 x QD16) : C2 & 1 & $C_8.C_2^2$ \\
 & 132 & (C2 x Q16) : C2 & 3 & $C_8.C_2^2, C_2 \times Q_{16}$ \\
 & 133 & ((C2 x Q8) : C2) : C2 & 1 & $C_8.C_2^2$ \\
 & 142 & (Q8 : C4) : C2 & 1 &  $C_8.C_2^2$ \\
 & 143 & C4 : Q16 & 3 & $C_8.C_2^2, C_2 \times Q_{16}$ \\
 & 145 & (C4 x Q8) : C2 & 1 & $C_8.C_2^2$ \\
 & 148 & (Q8 : C4) : C2 & 2 & $C_2 \times Q_{16}$ \\
 & 149 & ((C8 x C2) : C2) : C2 & 1 & $C_8.C_2^2$ \\
 & 151 & (Q8 : C4) : C2 & 2 & 2 copies of $C_2 \times Q_{16}$ \\
 & 154 & (C2 . ((C4 x C2) : C2) & 1 & \\
 & 155 & (C8 : C4) : C2 & 1 &  $C_8.C_2^2$ \\
 & 156 & Q8 : Q8 & 1 & $C_8.C_2^2$ \\
 & 158 & Q8 : Q8 & 2 & $C_2 \times Q_{16}$ \\
 & 160 & (C2 x C2) . (C2 x D8) & 1 & $C_8.C_2^2$ \\
 & 161 & ((C8 x C2) : C2) : C2 & 1 & $C_8.C_2^2$ \\
 & 164 & (Q8 : C4) : C2 & 1 & $C_8.C_2^2$ \\
 & 165 & (Q8 : C4) : C2 & 2 & $C_2 \times Q_{16}$ \\
 & 166 & (Q8 : C4) : C2 & 1 & $C_8.C_2^2$ \\
 & 168 & (C2 x C2) . (C2 x D8) & 2 &$C_2 \times Q_{16}$ \\
 & 175 & C4 : Q16 & 4 & 2 copies of $C_2 \times Q_{16}$ \\
 & 178 & (C4 : Q8) : C2 & 2 & $C_8.C_2^2$ \\
 & 181 & C8 : Q8 & 2 & $C_2 \times Q_{16}$ \\
 & 188 & $C_2 \times Q_{32}$ & 2 &\\
 & 191 & QD32 : C2 & 1 & \\
 & 252 & $C_2 \times C_2 \times Q_{16}$ & 4 & \\
 & 255 & C2 x ((C2 x Q8) : C2) & 2 & \\
 & 259 & (C2 x Q16) : C2 & 1 &  \\
\hline
\end{tabular} 
\caption{All groups $G$ with $s(G)>0$ of order $ n = 64$.} \label{groupsorder64}
\end{center}
\end{table}

\begin{table}
\begin{center}
\begin{tabular}{| c | c | c | c | c | }
\hline 
$n = 68$ & Index  &  Structure  & $s(G)$ & Quotients \\ 
\hline 
 & 1 & $\text{Dic}_{17}$ & 1 &  \\
\hline 
\hline
$n = 72$ & Index  &  Structure  & $s(G)$ & Quotients \\ 
\hline 
 & 4 & $\text{Dic}_{18}$ & 2 &  \\
 & 24 & (C3 x C3) : Q8 & 2 & 2 copies of $\text{Dic}_{6}$  \\
 & 26 & $C_3 \times \text{Dic}_6$ & 1 & \\
 & 31 & (C3 x C3) : Q8 & 4 & 4 copies of $\text{Dic}_{6}$  \\
\hline
\hline 
$n = 76$ & Index  &  Structure  & $s(G)$ & Quotients   \\ 
\hline 
\hline 
$n = 80$ & Index  &  Structure  & $s(G)$ & Quotients  \\ 
\hline 
 & 1 & C5 : C16 & 1 & $\text{Dic}_5$ \\
 & 8 & $\text{Dic}_{20}$ & 2 &\\
 & 9 & C2 x (C5 : C8) & 2 & \\
 & 10 & (C5 : C8) : C2 & 2 & 2 copies of $\text{Dic}_5$ \\
 & 11 & $C_4 \times \text{Dic}_{5}$ & 2 & \\
 & 12 & (C5 : C4) : C4 & 1 & $\text{Dic}_{10}$ \\
 & 13 & C20 : C4 & 3 & 2 copies of $\text{Dic}_5$ and 1 copy of $\text{Dic}_{10}$  \\
 & 18 & C5 : Q16 & 1 & $Q_{16}$  \\
 & 19 & (C10 x C2) : C4 & 2 & 2 copies of $\text{Dic}_5$  \\
 & 27 & $C_5 \times Q_{16}$ & 1 & \\
 & 33 & (C5 : C8) : C2 & 1 &  \\
 & 35 & $C_2 \times \text{Dic}_{10}$ & 2 & \\
 & 40 & (C4 x D10) : C2 & 1 &  \\
 & 43 & $C_2 \times C_2 \times \text{Dic}_{5}$ & 4 & \\
\hline 
\hline
$n = 84$ & Index  &  Structure  & $s(G)$ & Quotients  \\ 
\hline 
 & 5 & $\text{Dic}_{21}$ & 1 & \\
\hline 
\hline
$n = 88$ & Index  &  Structure  & $s(G)$ & Quotients  \\ 
\hline 
 & 3 & $\text{Dic}_{22}$ & 1 & \\
\hline 
\hline
$n = 92$ & Index  &  Structure  & $s(G)$ & Quotients   \\ 
\hline
\end{tabular}
\caption{All groups $H$ with $s(H)>0$ of order $68 \leq n \leq 92$.} \label{groupsorder65to95}
\end{center}
\end{table}

\begin{table}
\begin{center}
\begin{tabular}{| c | c | c | c | c | } 
\hline 
$n = 96$ & Index  &  Structure  & $s(G)$ & Quotients   \\ 
\hline 
 & 8 & $\text{Dic}_{24}$ & 2 & \\
 & 11 & C3 : (C4 : C8) & 1 & $\text{Dic}_6$ \\
 & 14 & C3 : (C8 : C4) & 2 & $Q_{16}, \text{Dic}_6$ \\
 & 15 & C3 : (C8 : C4) & 1 & $\text{Dic}_6$ \\
 & 17 & C3 : (Q8 : C4) & 1 & $Q_{16}$ \\
 & 21 & C3 : (C4 : C8) & 1 & $\text{Dic}_6$ \\
 & 23 & C3 : (Q8 : C4) & 2 & $Q_{16}, \text{Dic}_{12}$ \\
 & 24 & C3 : (C8 : C4) & 1 & $\text{Dic}_6$ \\
 & 25 & C3 : (C8 : C4) & 3 & $Q_{16}, \text{Dic}_6, \text{Dic}_{12}$ \\
 & 26 & C3 : (C4 . D8 = C4 . (C4 x C2)) & 1 & $\text{Dic}_6$ \\
 & 29 & C3 : (C4 . D8 = C4 . (C4 x C2)) & 1 & $\text{Dic}_6$ \\
 & 31 & C3 : (C2 . ((C4 x C2) : C2) & 1 & \\
 & 36 & C3 : Q32 & 1 & $Q_{32}$ \\
 & 38 & C3 : ((C4 x C2) : C4) & 1 & $\text{Dic}_6$ \\
 & 42 & C3 : (Q8 : C4) & 1 & $Q_{16}$ \\
 & 53 & C3 x (Q8 : C4) & 1 & $Q_{16}$ \\
 & 57 & C3 x (C8 : C4) & 1 & $Q_{16}$ \\
 & 63 & $C_3 \times Q_{32}$ & 1 & \\
 & 66 & SL(2,3) : C4 & 1 & $\tilde{O}$ \\
 & 75 & $C_4 \times \text{Dic}_6$ & 2 & \\
 & 76 & C3 : (C4 : Q8) & 4 & 4 copies of $\text{Dic}_6$ \\
 & 77 & C3 : ((C2 x C2) . (C2 x C2 x C2)) & 2 & 2 copies of $\text{Dic}_6$ \\
 & 85 & C3 : ((C2 x Q8) : C2) & 2 & 2 copies of $\text{Dic}_6$ \\
 & 95 & C3 : (C4 : Q8) & 2 & 2 copies of $\text{Dic}_6$ \\
 & 97 & C3 : ((C2 x C2) . (C2 x C2 x C2)) & 2 & 2 copies of $\text{Dic}_6$ \\
 & 112 & $C_2 \times \text{Dic}_{12}$ & 4 & \\
 & 116 & C3 : ((C2 x Q8) : C2) & 1 & $C_8 . C_2^2$ \\
 & 119 & C3 : ((C8 x C2) : C2) & 1 &  \\
 & 122 & C3 : ((C2 x Q8) : C2) & 1 & $C_8 . C_2^2$ \\
 & 124 & $Q_{16} \times S_3$ & 3 & \\
 & 130 & C2 x ((C3 : C4) : C4) & 2 & 2 copies of $\text{Dic}_6$ \\
 & 131 & C3 : ((C2 x Q8) : C2) & 2 & 2 copies of $\text{Dic}_6$ \\
 & 132 & C2 x (C12 : C4) & 2 & 2 copies of $\text{Dic}_6$ \\
 & 150 & C2 x (C3 : Q16) & 2 & 2 copies of $Q_{16}$ \\
 & 158 & C3 : ((C2 x Q8) : C2) & 1 & $C_8 . C_2^2$ \\
 & 181 & $C_6 \times Q_{16}$ & 2 & \\
 & 185 & A4 : Q8 & 1 & $\text{Dic}_6$ \\
 & 188 & $C_2 \times \tilde{O}$ & 2 & \\
 & 190 & (C2 x SL(2,3)) : C2 & 1 &  \\
 & 191 & (C2 . S4 = SL(2,3) . C2) : C2 & 2 & \\
 & 205 & $C_2 \times C_2 \times \text{Dic}_6$ & 4 & \\
 & 217 & C3 : ((C2 x Q8) : C2) & 1 &  \\
\hline 
$n = 100$ & Index  &  Structure  & $s(G)$ & Quotients   \\ 
\hline 
 & 1 & $\text{Dic}_{25}$ & 2 &\\
 & 6 & $C_5 \times \text{Dic}_{5}$ & 1 & \\
 & 7 & (C5 x C5) : C4 & 6 & 6 copies of $\text{Dic}_{5}$ \\
 & 10 & (C5 x C5) : C4 & 1 & $\text{Dic}_{5}$ \\
\hline
\end{tabular}
\caption{All groups $G$ with $s(G)>0$ of order $ 96 \leq n \leq 100$.} \label{groupsorder96to100}
\end{center}
\end{table}
\begin{table}
\begin{center}
\begin{tabular}{| c | c | c | c | c | c |}
\hline
$n$ & Index & Structure & s \\
\hline
16 & 9 & $Q_{16}$ & 1 \\
20 & 1 & $\text{Dic}_5$ & 1 \\
24 & 4 & $\text{Dic}_6$ & 1 \\
32 & 20 & $Q_{32}$ & 1 \\
32 & 44 & $C_8.C_2^2$ & 1 \\
40 & 4 & $\text{Dic}_{10}$ & 1 \\
48 & 28 & $\tilde{O}$ & 1 \\
52 & 1 & $\text{Dic}_{13}$ & 1 \\
56 & 3 & $\text{Dic}_{14}$ & 1 \\
64 & 54 & $Q_{64}$ & 1 \\
64 & 154 & (C2 . ((C4 x C2) : C2) = (C2 x C2) . (C4 x C2)) : C2 & 1 \\
64 & 191 & QD32 : C2 & 1 \\
64 & 259 & (C2 x Q16) : C2 & 1 \\
68 & 1 & $\text{Dic}_{17}$ & 1 \\
80 & 33 & $C_2^2.F_5$ & 1 \\
80 & 40 & (C4 x D10) : C2 & 1 \\
84 & 5 & $\text{Dic}_{21}$ & 1 \\
88 & 3 & $\text{Dic}_{22}$ & 1 \\
96 & 31 & C3 : (C2 . ((C4 x C2) : C2) = (C2 x C2) . (C4 x C2)) & 1 \\
96 & 119 & C3 : ((C8 x C2) : C2) & 1 \\
96 & 190 & (C2 x SL(2,3)) : C2 & 1 \\
96 & 191 & (C2 . S4 = SL(2,3) . C2) : C2 & 2 \\
96 & 217 & C3 : ((C2 x Q8) : C2) & 1 \\
104 & 4 & $\text{Dic}_{26}$ & 1 \\
116 & 1 & $\text{Dic}_{29}$ & 1 \\
120 & 5 & $\tilde{I}$ & 1 \\
120 & 8 & (C3 : C4) x D10 & 1 \\
128 & 66 & Q16 : C8 & 1 \\
128 & 72 & (C2 x Q16) : C4 & 1 \\
128 & 74 & C2 . (((C8 x C2) : C2) : C2) = (C8 x C2) . (C4 x C2) & 1 \\
128 & 82 & (C2 x C2) . ((C8 x C2) : C2) = (C8 x C2) . (C4 x C2) & 1 \\
128 & 90 & C2 . ((C4 : C8) : C2) = (C4 x C4) . (C4 x C2) & 2 \\
128 & 137 & ((C2 x C2) . ((C4 x C2) : C2) = (C4 x C2) . (C4 x C2)) : C2 & 1 \\
128 & 143 & C2 . (((C2 x C2 x C2) : C4) : C2) = (C4 x C4) . (C4 x C2) & 2 \\
128 & 152 & C2 . ((C16 x C2) : C2) = C16 . (C4 x C2) & 1 \\
128 & 163 & $Q_{128}$ & 1 \\
128 & 634 & ((C8 x C2) : C4) : C2 & 1 \\
128 & 637 & ((Q8 : C4) : C2) : C2 & 1 \\
128 & 879 & C2 x (C2 . ((C8 x C2) : C2) = C8 . (C4 x C2)) & 2 \\
128 & 912 & (Q16 : C4) : C2 & 1 \\
128 & 927 & ((C2 x Q16) : C2) : C2 & 1 \\
128 & 946 & (C2 . ((C8 x C2) : C2) = C8 . (C4 x C2)) : C2 & 2 \\
128 & 954 & (C2 . ((C8 x C2) : C2) = C8 . (C4 x C2)) : C2 & 1 \\
128 & 971 & (C2 . ((C8 x C2) : C2) = C8 . (C4 x C2)) : C2 & 1 \\
\hline
\end{tabular} 
\caption{All groups $G$ with $s(G)>0$ of order $ n \leq 180$ with no quotient with $s > 0$, part 1.} \label{groupswithnoquotients128}
\end{center}
\end{table}

\clearpage

\begin{table}
\begin{center}
\begin{tabular}{| c | c | c | c | c | c |}
\hline
$n$ & Index & Structure & s \\
\hline
128 & 996 & QD64 : C2 & 1 \\
128 & 2025 & (((C2 x Q8) : C2) : C2) : C2 & 2 \\
128 & 2149 & (C2 x Q32) : C2 & 1 \\
128 & 2318 & (C2 x ((C2 x Q8) : C2)) : C2 & 1 \\
132 & 3 & $\text{Dic}_{33}$ & 1 \\
136 & 3 & C17 : C8 & 1 \\
136 & 4 & $\text{Dic}_{34}$ & 1 \\
144 & 15 & C9 : QD16 & 1 \\
144 & 43 & Q8 x D18 & 1 \\
144 & 59 & (C3 x C3) : QD16 & 1 \\
144 & 118 & (C3 x C3) : QD16 & 1 \\
144 & 138 & (C3 x C3) : ((C4 x C2) : C2) & 1 \\
148 & 1 & $\text{Dic}_{37}$  & 1 \\
152 & 3 & $\text{Dic}_{38}$ & 1 \\
156 & 1 & C13 : C12 & 1 \\
160 & 31 & C5 : (C2 . ((C4 x C2) : C2) = (C2 x C2) . (C4 x C2)) & 1 \\
160 & 80 & C5 : (C2 . ((C4 x C2) : C2) = (C2 x C2) . (C4 x C2)) & 1 \\
160 & 133 & C5 : ((C8 x C2) : C2) & 1 \\
160 & 225 & C5 : ((C2 x Q8) : C2) & 1 \\
164 & 1 & $\text{Dic}_{41}$ & 1 \\
168 & 7 & C7 : (C3 x Q8) & 1 \\
168 & 12 & (C3 : C4) x D14 & 1 \\
168 & 15 & C3 : ((C14 x C2) : C2) & 1 \\
\hline
\end{tabular} 
\caption{All groups $G$ with $s(G)>0$ of order $n \leq 180$ with no quotient with $s > 0$, part 2.} \label{groupswithnoquotients180}
\end{center}
\end{table}

\subsection{Methodology}
All computations have been done on a Microsoft Surface Pro 6 with Windows 10 as operating system. The software used was GAP 4.11.0 \cite{GAP4}.

For the computation of groups with non-trivial $s$ the following code was used.

\begin{lstlisting}[language=GAP]
output := OutputTextFile( "output.txt" , true );

SearchAllGroupsForS := function(n)
local j, s, G;
for j in [1..Length(AllSmallGroups(n))] do
G := SmallGroup(n,j);
s := sOfGroup(G);
if s>0 then AppendTo(output, " & ",j," & ",StructureDescription(G)," & ",s,
      " & \\","\\","\n"); fi;
od;
end;

for k in Filtered([4..100], i -> GcdInt(i,4)=4) do
# k is the order of the groups searched. 
# Note that we only have s>0 for groups of order divisible by 4

AppendTo(output, "\\","hline \n","$n = ", k,"$ & Index  &  Structure  & $s(G)$ 
      & Quotients \\","\\ \n","\\","hline \n");
SearchAllGroupsForS(k);
od;
\end{lstlisting}

For the computation of those groups with non-trivial $s$ without a smaller quotient with $s>0$ the following code was used. Here \texttt{ElementGrps} is a list of the \texttt{SmallGroupLibrary} IDs of those groups with non-trivial $s$ without quotients with that property. The algorithm presentated checks all quotients of a given group with $s>0$ against smaller groups that already appear in \texttt{ElementGrps}. This is done by computing \texttt{NormalSubgroups(G)} and for each normal subgroup \texttt{N} compute the ID of the quotient \texttt{G/N}, then compare this new ID with all IDs already in \texttt{ElementGrps}. If no match is found, add the ID of \texttt{G} to \texttt{ElementGrps}.

\begin{lstlisting}[language=GAP]
ElementGrps := []

SearchAllBuildingBlocksForS := function(n)
local j, s, G;
for j in [1..Length(AllSmallGroups(n))] do
G := SmallGroup(n,j);
s := sOfGroup(G);
if s>0 then 
  if not ForAny( NormalSubgroups(G), 
     N ->  IdSmallGroup( FactorGroup( G, N ) ) in ElementGrps )
     then Add(ElementGrps, [n,j]); fi;
  fi;          
od;
end;

for k in Filtered([4..180], i -> GcdInt(i,4)=4) do
# k is the order of the groups searched. 
# Note that we only have s>0 for groups of order divisible by 4
SearchAllBuildingBlocksForS(k);
od;

Print( ElementGrps );
\end{lstlisting}

We want to remark that the computations for orders above 100, especially at order 128, tended to be quite slow. The actual computation has been done in steps for fixed intervals of orders rather than all at once.

\begingroup
\setlength{\emergencystretch}{8em}
\printbibliography
\endgroup

\end{document}